\documentclass[11pt,a4paper]{article}
\usepackage{amsmath,amssymb,amsthm,mathtools}
\usepackage[a4paper,margin=2.7cm]{geometry}
\usepackage{microtype}
\usepackage{enumitem}
\usepackage{hyperref}
\hypersetup{colorlinks=true,linkcolor=black,citecolor=black,urlcolor=black}

\newtheorem{theorem}{Theorem}[section]
\newtheorem{proposition}[theorem]{Proposition}
\newtheorem{lemma}[theorem]{Lemma}
\newtheorem{remark}[theorem]{Remark}

\title{Two Integral Sliding-Window Inequalities for Maximal Convolutions}
\author{Gangsong Leng\thanks{East China Normal University, email: \texttt{lenggangsong@163.com}}\quad Cheng Li\thanks{School of Mathematical Sciences, Fudan University, email: \texttt{24300180118@m.fudan.edu.cn}}}
\date{}

\begin{document}
\maketitle

\begin{abstract}
We prove two sliding-window inequalities for maximal convolutions. The first concerns the multiplicative maximal convolution. If $f$ and $g$ are nonnegative continuous functions on $[0,A]$ and $[0,B]$, respectively, define
\[
    h(x)=\max_{\substack{0\le u\le A\\0\le x-u\le B}} f(u)g(x-u),\qquad 0\le x\le A+B.
\]
Then there exists a window $[a,a+B]$ of length $B$ such that
\[
    \frac1B\int_a^{a+B}h(x)\,dx\ge
    \left(\frac1A\int_0^A f(x)\,dx\right)
    \left(\frac1B\int_0^B g(x)\,dx\right).
\]
The second concerns the additive maximal convolution. Let $f$ and $g$ be nonnegative continuous functions on $[0,C]$, and define
\[
    H(x)=\max_{\substack{0\le u\le C\\0\le x-u\le C}}\{f(u)+g(x-u)\},\qquad 0\le x\le 2C.
\]
Then, for every $p\ge1$, there exists a window $[a,a+C]$ of length $C$ such that
\[
    \left(\int_a^{a+C}H(x)^p\,dx\right)^{1/p}
    \ge
    \left(\int_0^C f(x)^p\,dx\right)^{1/p}
    +
    \left(\int_0^C g(x)^p\,dx\right)^{1/p}.
\]
We also record discrete analogues. The main point is that, in one-dimensional maximal-convolution settings, certain global Brunn--Minkowski or Pr\'ekopa--Leindler type phenomena admit natural sliding-window localizations.
\end{abstract}



\section{Background and main results}

Ordinary convolution sums all contributions with $u+v=x$. In contrast, maximal convolution keeps only the largest contribution. We consider two versions. The first is multiplicative:
\[
    h(x)=\max_{u+v=x}f(u)g(v),
\]
and the second is additive:
\[
    H(x)=\max_{u+v=x}\{f(u)+g(v)\}.
\]
These operations are functional analogues of sumsets. If $f$ and $g$ are indicator functions of two sets, then the multiplicative maximal convolution is the indicator function of their sumset. The additive maximal convolution is the corresponding max-plus operation after taking logarithms in the multiplicative setting.

Maximal convolutions are naturally related to the Brunn--Minkowski, Pr\'ekopa--Leindler and Borell--Brascamp--Lieb circles of ideas. For the Brunn--Minkowski inequality and its many equivalent forms and applications, we refer to Gardner's survey \cite{Gardner2002}. For the Pr\'ekopa--Leindler inequality and Brascamp--Lieb type extensions, see Brascamp and Lieb \cite{BrascampLieb1976}. In the discrete direction, the Ahlswede--Daykin four functions theorem \cite{AhlswedeDaykin1978} is an important discrete functional inequality and is often viewed as one representative form of the discrete Pr\'ekopa--Leindler philosophy.

The results below, however, are not direct special cases of the standard Pr\'ekopa--Leindler inequality. The classical Pr\'ekopa--Leindler theorem concerns control at weighted midpoints $(1-\lambda)x+\lambda y$ through geometric means. Here we deal instead with maximal convolutions generated by the sum constraint $u+v=x$, and we prove that a fixed-length sliding window already carries the required lower bound. Thus the connection with Pr\'ekopa--Leindler type inequalities should be understood as one of common background rather than as a direct implication.

The main novelty is the sliding-window localization. A usual global maximal-convolution, Brunn--Minkowski or Pr\'ekopa--Leindler type statement controls a quantity on the whole convolution interval. In contrast, our additive theorem shows that a single interval of the original length already carries the full lower bound $\|f\|_p+\|g\|_p$. The multiplicative theorem gives a parallel local average statement for the multiplicative maximal convolution. These assertions should not be regarded as consequences of standard Pr\'ekopa--Leindler theory.

It is useful to stress why the additive theorem is not merely a consequence of averaging over all windows. A direct window-averaging argument gives only the weaker estimate
\[
    \left(\int_a^{a+C} H(x)^p\,dx\right)^{1/p}
    \ge
    \left(\|f\|_p^p+\|g\|_p^p\right)^{1/p}
\]
for some window $[a,a+C]$. It does not yield the stronger lower bound
\[
    \|f\|_p+\|g\|_p.
\]
The proof of Theorem \ref{thm:add-int} uses instead a maximum-window argument together with a secant-line lemma. This is precisely where the sliding-window phenomenon goes beyond elementary window averaging.

Our first result is an integral sliding-window inequality for the multiplicative maximal convolution.

\begin{theorem}[Integral sliding-window inequality for multiplicative maximal convolution]\label{thm:mult-int}
Let $A,B>0$. Let $f$ and $g$ be nonnegative continuous functions on $[0,A]$ and $[0,B]$, respectively. For $0\le x\le A+B$, define
\[
    h(x)=\max_{\substack{0\le u\le A\\0\le x-u\le B}} f(u)g(x-u).
\]
For each $x\in[0,A+B]$, the feasible set is nonempty and compact, so the maximum exists. Then there exists $a\in[0,A]$ such that
\[
    \int_a^{a+B}h(x)\,dx
    \ge
    \frac1A
    \left(\int_0^A f(x)\,dx\right)
    \left(\int_0^B g(x)\,dx\right).
\]
Equivalently,
\[
    \frac1B\int_a^{a+B}h(x)\,dx
    \ge
    \left(\frac1A\int_0^A f(x)\,dx\right)
    \left(\frac1B\int_0^B g(x)\,dx\right).
\]
\end{theorem}

The proof is a continuous double-counting argument. Integrating the window integrals over all window origins turns the problem into an integral over the rectangle $[0,A]\times[0,B]$.

The second result is the main theorem of the paper. It gives an integral sliding-window inequality for the additive maximal convolution. Its discrete form was first obtained by the second-named author in a short note on a reverse Minkowski inequality \cite{Li2024}.

\begin{theorem}[Integral sliding-window inequality for additive maximal convolution]\label{thm:add-int}
Let $C>0$, $p\ge1$, and let $f,g$ be nonnegative continuous functions on $[0,C]$. For $0\le x\le2C$, define
\[
    H(x)=\max_{\substack{0\le u\le C\\0\le x-u\le C}}\{f(u)+g(x-u)\}.
\]
For each $x\in[0,2C]$, the feasible set
\[
    \{u:0\le u\le C,\ 0\le x-u\le C\}
\]
is nonempty and compact, so the maximum exists. Then there exists $a\in[0,C]$ such that
\[
    \left(\int_a^{a+C}H(x)^p\,dx\right)^{1/p}
    \ge
    \left(\int_0^C f(x)^p\,dx\right)^{1/p}
    +
    \left(\int_0^C g(x)^p\,dx\right)^{1/p}.
\]
\end{theorem}

If the integral on the left is taken over the whole interval $[0,2C]$, then Theorem \ref{thm:add-int} immediately yields a global $L^p$ inequality for additive maximal convolution. For $p=2$, this global inequality has a direct interpretation through the two-dimensional Brunn--Minkowski inequality; for positive integer $p$, it can similarly be interpreted through the $p$-dimensional Brunn--Minkowski inequality. Theorem \ref{thm:add-int} is more precise: it shows that one local window of length $C$ inside the global interval already supports the sum of the two $L^p$ norms on the right-hand side.

We also record two discrete counterparts. They should be viewed as the corresponding models on integer grids with counting measure. Since the passage between continuous intervals and integer grids involves minor boundary differences, we state them as propositions and give direct discrete proofs in the next section rather than deriving them by a discretization argument.

\begin{proposition}[Discrete window form for multiplicative maximal convolution]\label{prop:mult-discrete}
Let
\[
    a_0,\ldots,a_m\ge0,
    \qquad
    b_0,\ldots,b_n\ge0,
\]
and define
\[
    c_k=\max_{\substack{i+j=k\\0\le i\le m,\ 0\le j\le n}}a_i b_j,
    \qquad k=0,1,\ldots,m+n.
\]
Then there exists $r\in\{0,1,\ldots,m\}$ such that
\[
    \frac1{n+1}\sum_{k=r}^{r+n}c_k
    \ge
    \left(\frac1{m+1}\sum_{i=0}^m a_i\right)
    \left(\frac1{n+1}\sum_{j=0}^n b_j\right).
\]
Symmetrically, there exists $s\in\{0,1,\ldots,n\}$ such that
\[
    \frac1{m+1}\sum_{k=s}^{s+m}c_k
    \ge
    \left(\frac1{m+1}\sum_{i=0}^m a_i\right)
    \left(\frac1{n+1}\sum_{j=0}^n b_j\right).
\]
\end{proposition}

\begin{proposition}[Discrete window form for additive maximal convolution]\label{prop:add-discrete}
Let $p\ge1$, and let $s_0,\ldots,s_n,t_0,\ldots,t_n$ be nonnegative real numbers. Define
\[
    r_k=\max_{i+j=k}(s_i+t_j),
    \qquad k=0,1,\ldots,2n.
\]
Then there exists $m\in\{0,1,\ldots,n\}$ such that
\[
    \left(\sum_{k=m}^{m+n}r_k^p\right)^{1/p}
    \ge
    \left(\sum_{i=0}^n s_i^p\right)^{1/p}
    +
    \left(\sum_{j=0}^n t_j^p\right)^{1/p}.
\]
More generally, the same conclusion holds under the weaker assumption
\[
    r_{i+j}\ge s_i+t_j\qquad (0\le i,j\le n).
\]
\end{proposition}

\begin{remark}
The two discrete corollaries should be compared with global inequalities in different ways. The multiplicative discrete window form is a local average statement; it has a different emphasis from the corresponding global average inequality, so there is no simple ordering by strength. The additive discrete window form, on the other hand, is a genuine localization of the global additive maximal-convolution $\ell^p$ inequality, since a single window already gives the same lower bound as the global form.
\end{remark}

\section{Proofs}

We first prove the multiplicative integral theorem.

\begin{proof}[Proof of Theorem \ref{thm:mult-int}]
For $a\in[0,A]$, set
\[
    W(a)=\int_a^{a+B}h(x)\,dx.
\]
The function $h$ is continuous by the maximum theorem, hence $W$ is continuous as well, although only integrability is needed below.
Then
\[
\begin{aligned}
    \int_0^A W(a)\,da
    &=\int_0^A\int_a^{a+B}h(x)\,dx\,da.
\end{aligned}
\]
Make the change of variables
\[
    u=a,
    \qquad
    v=x-a.
\]
When $0\le a\le A$ and $a\le x\le a+B$, this is exactly the region
\[
    0\le u\le A,
    \qquad
    0\le v\le B,
    \qquad
    x=u+v.
\]
The Jacobian is equal to $1$. Hence
\[
    \int_0^A\int_a^{a+B}h(x)\,dx\,da
    =
    \int_0^A\int_0^B h(u+v)\,dv\,du.
\]
By the definition of $h$,
\[
    h(u+v)\ge f(u)g(v).
\]
Therefore
\[
\begin{aligned}
    \int_0^A W(a)\,da
    &\ge \int_0^A\int_0^B f(u)g(v)\,dv\,du  \\
    &=\left(\int_0^A f(u)\,du\right)
      \left(\int_0^B g(v)\,dv\right).
\end{aligned}
\]
Thus the average of $W(a)$ over $[0,A]$ is at least
\[
    \frac1A
    \left(\int_0^A f\right)
    \left(\int_0^B g\right).
\]
Consequently there exists $a\in[0,A]$ such that
\[
    W(a)\ge
    \frac1A
    \left(\int_0^A f\right)
    \left(\int_0^B g\right),
\]
as required.
\end{proof}

We next prove the additive integral theorem. We begin with a simple secant-line lemma.

\begin{lemma}[Secant-line lemma]\label{lem:secant}
Let $x,y>0$, $0\le z\le y$, and $p\ge1$. Then
\[
    (x+z)^p
    \ge
    x^p+\frac{(x+y)^p-x^p}{y^p}z^p.
\]
\end{lemma}

\begin{proof}
For $p=1$ the assertion is an equality. Assume $p>1$. For $z>0$, define
\[
    \Phi(z)=\frac{(x+z)^p-x^p}{z^p}.
\]
Writing $u=z/x$, we have
\[
    \Phi(z)=\frac{(1+u)^p-1}{u^p}.
\]
A direct differentiation gives
\[
    \frac{d}{du}\left(\frac{(1+u)^p-1}{u^p}\right)
    =p u^{-p-1}\left(1-(1+u)^{p-1}\right)\le0.
\]
Hence $\Phi$ is decreasing. Since $0<z\le y$, it follows that
\[
    \frac{(x+z)^p-x^p}{z^p}
    \ge
    \frac{(x+y)^p-x^p}{y^p},
\]
which is the desired inequality. The case $z=0$ is immediate.
\end{proof}

\begin{proof}[Proof of Theorem \ref{thm:add-int}]
Set
\[
    S=\max_{0\le x\le C}f(x),
    \qquad
    T=\max_{0\le x\le C}g(x).
\]
If $S=0$ or $T=0$, the conclusion is immediate. For example, if $S=0$, then $f\equiv0$. For any $a\in[0,C]$ and $x\in[a,a+C]$, write $x=a+v$ with $v\in[0,C]$. Then
\[
    H(a+v)\ge f(a)+g(v)=g(v),
\]
and hence
\[
    \left(\int_a^{a+C}H(x)^p\,dx\right)^{1/p}
    \ge \left(\int_0^C g(v)^p\,dv\right)^{1/p}.
\]
The case $T=0$ is analogous. We may therefore assume that $S,T>0$.

Choose $x_0,y_0\in[0,C]$ such that
\[
    f(x_0)=S,
    \qquad
    g(y_0)=T.
\]
For $a\in[0,C]$, define
\[
    W(a)=\left(\int_a^{a+C}H(x)^p\,dx\right)^{1/p}.
\]
Since $f$ and $g$ are continuous on compact intervals, the maximum theorem implies that $H$ is continuous. Hence $W$ is continuous, and there exists $a_0\in[0,C]$ such that
\[
    W(a_0)=\max_{0\le a\le C}W(a).
\]

We first use the point where $f$ attains its maximum. For every $0\le v\le C$,
\[
    H(x_0+v)\ge f(x_0)+g(v)=S+g(v).
\]
Thus
\[
    W(x_0)^p
    \ge \int_0^C (S+g(v))^p\,dv.
\]
Since $W(a_0)\ge W(x_0)$, we get
\[
    W(a_0)^p
    \ge \int_0^C (S+g(v))^p\,dv.       \tag{2.1}
\]
Similarly, using the point where $g$ attains its maximum, we obtain
\[
    W(a_0)^p
    \ge \int_0^C (T+f(u))^p\,du.       \tag{2.2}
\]

Let
\[
    \left(\int_0^C f(u)^p\,du\right)^{1/p}=\lambda S,
    \qquad
    \left(\int_0^C g(v)^p\,dv\right)^{1/p}=\mu T.
\]
Since $S$ and $T$ are maximum values, we have
\[
    0\le \lambda^p,\mu^p\le C.
\]
By Lemma \ref{lem:secant}, applied with $x=S$, $y=T$ and $z=g(v)$, we have for every $v\in[0,C]$,
\[
    (S+g(v))^p
    \ge
    S^p+\frac{(S+T)^p-S^p}{T^p}g(v)^p.
\]
Integrating this inequality yields
\[
\begin{aligned}
    \int_0^C(S+g(v))^p\,dv
    &\ge C S^p+\frac{(S+T)^p-S^p}{T^p}\int_0^C g(v)^p\,dv \\
    &=C S^p+\mu^p\bigl((S+T)^p-S^p\bigr) \\
    &=\mu^p(S+T)^p+(C-\mu^p)S^p \\
    &\ge \mu^p(S+T)^p .
\end{aligned}
\]
Together with (2.1), this gives
\[
    W(a_0)\ge \mu(S+T).       \tag{2.3}
\]
In the same way, (2.2) gives
\[
    W(a_0)\ge \lambda(S+T).   \tag{2.4}
\]
Hence
\[
    W(a_0)\ge \max\{\lambda,\mu\}(S+T).
\]
Finally,
\[
    \max\{\lambda,\mu\}(S+T)
    \ge \lambda S+\mu T.
\]
Consequently,
\[
\begin{aligned}
    W(a_0)
    &\ge \lambda S+\mu T \\
    &=\left(\int_0^C f(u)^p\,du\right)^{1/p}
      +\left(\int_0^C g(v)^p\,dv\right)^{1/p}.
\end{aligned}
\]
This proves the theorem.
\end{proof}

It remains to prove the two discrete propositions.

\begin{proof}[Proof of Proposition \ref{prop:mult-discrete}]
Set
\[
    W_r=\sum_{k=r}^{r+n}c_k,
    \qquad r=0,1,\ldots,m.
\]
Then
\[
    \sum_{r=0}^m W_r
    =\sum_{r=0}^m\sum_{k=r}^{r+n}c_k.
\]
For fixed $k$, the number of windows in which $c_k$ occurs is precisely the number of pairs $(i,j)$ satisfying
\[
    i+j=k,
    \qquad 0\le i\le m,
    \qquad 0\le j\le n.
\]
Therefore
\[
    \sum_{r=0}^m W_r
    =\sum_{i=0}^m\sum_{j=0}^n c_{i+j}.
\]
Since $c_{i+j}\ge a_i b_j$, we have
\[
    \sum_{r=0}^m W_r
    \ge
    \sum_{i=0}^m\sum_{j=0}^n a_i b_j
    =
    \left(\sum_{i=0}^m a_i\right)
    \left(\sum_{j=0}^n b_j\right).
\]
Hence at least one $r$ satisfies
\[
    W_r\ge
    \frac1{m+1}
    \left(\sum_{i=0}^m a_i\right)
    \left(\sum_{j=0}^n b_j\right).
\]
Dividing by $n+1$ gives the first assertion. The symmetric assertion is proved in the same way.
\end{proof}

\begin{proof}[Proof of Proposition \ref{prop:add-discrete}]
We prove the more general form under the assumption $r_{i+j}\ge s_i+t_j$. If one of the two sequences is identically zero, the result is immediate. Assume
\[
    S=\max_i s_i>0,
    \qquad
    T=\max_j t_j>0.
\]
Choose indices $\alpha$ and $\beta$ such that
\[
    s_\alpha=S,
    \qquad
    t_\beta=T.
\]
Define
\[
    L_m=\left(\sum_{k=m}^{m+n}r_k^p\right)^{1/p},
    \qquad m=0,1,\ldots,n,
\]
and choose $m_0$ so that $L_{m_0}=\max_m L_m$. Let
\[
    \left(\sum_{i=0}^n s_i^p\right)^{1/p}=\lambda S,
    \qquad
    \left(\sum_{j=0}^n t_j^p\right)^{1/p}=\mu T.
\]
Since $S$ and $T$ are maximum entries, $1\le \lambda,\mu\le (n+1)^{1/p}$.

From $r_{\alpha+\ell}\ge S+t_\ell$, we have
\[
    L_\alpha^p\ge\sum_{\ell=0}^n(S+t_\ell)^p.
\]
By Lemma \ref{lem:secant},
\[
    (S+t_\ell)^p
    \ge
    S^p+\frac{(S+T)^p-S^p}{T^p}t_\ell^p.
\]
Summing over $\ell$ gives
\[
\begin{aligned}
    \sum_{\ell=0}^n(S+t_\ell)^p
    &\ge (n+1)S^p+
      \frac{(S+T)^p-S^p}{T^p}
      \sum_{\ell=0}^n t_\ell^p \\
    &= (n+1)S^p+\mu^p\bigl((S+T)^p-S^p\bigr) \\
    &= \mu^p(S+T)^p+(n+1-\mu^p)S^p \\
    &\ge \mu^p(S+T)^p .
\end{aligned}
\]
Thus
\[
    L_{m_0}\ge L_\alpha\ge \mu(S+T).
\]
Similarly, from $r_{\beta+\ell}\ge s_\ell+T$, we get
\[
    L_{m_0}\ge L_\beta\ge \lambda(S+T).
\]
Therefore
\[
    L_{m_0}\ge \max\{\lambda,\mu\}(S+T)
    \ge \lambda S+\mu T.
\]
The right-hand side is exactly
\[
    \left(\sum_{i=0}^n s_i^p\right)^{1/p}
    +
    \left(\sum_{j=0}^n t_j^p\right)^{1/p}.
\]
The proof is complete.
\end{proof}


\begin{thebibliography}{9}

\bibitem{Li2024}
C. Li, A reverse Minkowski inequality, \emph{NSMath Student Column}, 2024, revised August 1, 2024, \url{http://www.nsmath.cn/xszl} (in Chinese).

\bibitem{AhlswedeDaykin1978}
R. Ahlswede and D. E. Daykin, An inequality for the weights of two families of sets, their unions and intersections, \emph{Zeitschrift f\"ur Wahrscheinlichkeitstheorie und verwandte Gebiete} \textbf{43} (1978), 183--185.

\bibitem{Borell1975}
C. Borell, Convex set functions in $d$-space, \emph{Periodica Mathematica Hungarica} \textbf{6} (1975), no. 2, 111--136.

\bibitem{BrascampLieb1976}
H. J. Brascamp and E. H. Lieb, On extensions of the Brunn--Minkowski and Pr\'ekopa--Leindler theorems, including inequalities for log concave functions, and with an application to the diffusion equation, \emph{Journal of Functional Analysis} \textbf{22} (1976), 366--389.

\bibitem{Gardner2002}
R. J. Gardner, The Brunn--Minkowski inequality, \emph{Bulletin of the American Mathematical Society} \textbf{39} (2002), 355--405.

\end{thebibliography}
\end{document}